\documentclass[oneside,reqno,12pt]{amsart}
\usepackage[greek,english]{babel}
\usepackage{amsthm}
\usepackage{amsbsy}
\usepackage{amsfonts}
\usepackage{graphicx}
\usepackage{hyperref}
\usepackage{appendix}
\hypersetup{colorlinks=true, citecolor=blue, linkcolor=red}

\newtheorem{thm}{Theorem}

\newtheorem{lem}{Lemma}

\newtheorem{rem}{Remark}

\begin{document}
\title[A  Liouville theorem for  ancient solutions]{A Liouville theorem for  ancient solutions to a  semilinear  heat equation and its elliptic counterpart}
\author{Christos Sourdis}
\address{National and Kapodistrian University of Athens, Department of Mathematics, Greece.
}
\email{sourdis@uoc.gr}
\date{\today}
\begin{abstract}
 We establish the nonexistence of nontrivial ancient solutions
to the nonlinear heat equation
$u_t=\Delta u+|u|^{p-1}u$ which are smaller in absolute value than the self-similar radial singular steady state, provided that the exponent $p$ is strictly between Serrin's exponent  and that of Joseph and Lundgren.
This result was previously established by Fila and Yanagida [Tohoku Math. J.  (2011)] by  using forward self-similar solutions as barriers.   In contrast, we apply a sweeping argument with a family of time independent weak supersolutions. Our approach naturally lends itself to yield an analogous Liouville type result for the steady state problem in higher dimensions. In fact, in the case of the critical Sobolev exponent we show the validity of our  results for solutions that are smaller in absolute value than
a 'Delaunay'-type singular solution.\end{abstract}
 \maketitle

\section{Introduction}

We consider classical solutions to the semilinear equation
\begin{equation}\label{eqAncient}
        u_t=\Delta u+|u|^{p-1}u,\ \ x\in\mathbb{R}^N,\ t\leq 0,
\end{equation}
with $p>1$. For obvious reasons, such solutions are frequently called ancient.
Our interest will be in conditions which imply $u\equiv 0$,  a Liouville type theorem that is.

In the  past few decades there have been intensive studies of Liouville type theorems for the equation in (\ref{eqAncient}), either when $t\leq 0$, $t\in \mathbb{R}$ (entire solutions) or $t\geq 0$ (global solutions).
At the
same time, these have emerged as a fundamental tool in deriving various
qualitative properties of the solutions to the corresponding Cauchy problem in a general domain or for a nonlinearity that behaves like a power as $u\to \infty$. The best general reference here is the monograph \cite{qs}. For a recent account of the theory and some further developments, we refer to  \cite{polacikQui}.

The following three exponents play an important role in the study of the  equation in (\ref{eqAncient}):
\[\textrm{Serrin's exponent}\ \ p_{sg}=\frac{N}{N-2}\ \ \textrm{if}\ N\geq 3,\ \ p_{sg}=\infty \ \textrm{if}\ N=1,2;\]
\[\textrm{the critical Sobolev exponent}\ \ p_S=\frac{N+2}{N-2}\ \ \textrm{if}\ N\geq 3,\ \ p_S=\infty \ \textrm{if}\ N=1,2;\]
 \[
\textrm{the Joseph-Lundgren exponent}\ \  p_{JL}=\left\{\begin{array}{ll}
              \frac{(N-2)^2-4N+8\sqrt{N-1}}{(N-2)(N-10)}&\textrm{if}\ N> 10, \\
      \infty &\textrm{if}\ N\leq 10.
            \end{array}
    \right.
\]
We note that $p_{sg}<p_S<p_{JL}$ if $N\geq 3$.  These exponents arise naturally in the study of the ordinary differential equation that is satisfied by the positive radial steady states \cite{joseph,qs}. In this regard, let us list some well known properties   which we will need in the sequel.
First, for $p>p_{sg}$  there exists an  explicit radial singular steady state
\begin{equation}\label{eqSingular}
  \varphi_\infty(x)=L|x|^{-2/(p-1)}\ \ \textrm{with}\ \ L=\left(\frac{2}{p-1}\left(N-2-\frac{2}{p-1}\right) \right)^{1/(p-1)}.
\end{equation}
For any $p>1$, the radial ODE for the steady states admits a unique solution $\Phi$ such that $\Phi(0)=1$, $\Phi_r(0)=0$. This solution is defined in a maximal interval
of the form $[0,R_{max})$ with $0<R_{max}\leq \infty$ and is decreasing as long as it stays positive.
The following qualitative properties of $\Phi$ will be useful for our purposes.
\begin{itemize}
\item If $p\in (1,p_S)$, then $\Phi$ has a first root $\rho>0$. Actually, there are no positive steady states in this regime (see \cite{GidiGracie}). In fact, under the further restriction that $p\in (p_{sg},p_S)$, it intersects twice with $\varphi_{\infty}$ in $(0,\rho)$ (see in particular \cite[Fig. 2]{joseph}). We point out that these intersections are transverse thanks to the uniqueness of solutions to the corresponding IVP;
  \item If $p=p_S$, then $R_{max}=\infty$, $\Phi>0$, $\Phi(\infty)=0$ and $\Phi$ has exactly two  intersections with $\varphi_\infty$. We point out that this $\Phi$ has a  simple explicit formula and decays to zero faster than $\varphi_\infty$ as $r\to \infty$ (see also the discussion following (\ref{eqODEAI}) below);
  \item If $p\in (p_S,p_{JL})$, then $R_{max}=\infty$, $\Phi>0$, $\Phi(\infty)=0$ and $\Phi$ has infinitely many   intersections with $\varphi_\infty$. Moreover, $\Phi/\varphi_\infty\to 1$ as $r\to \infty$;
  \item If $p\geq p_{JL}$, then $R_{max}=\infty$, $\Phi>0$, $\Phi(\infty)=0$ and $\Phi<\varphi_\infty$. Moreover, $\Phi/\varphi_\infty\to 1$ as $r\to \infty$ (detailed information on this asymptotic behaviour can be found in \cite{GracieNW}).
\end{itemize}
 We emphasize that the rescaling
\begin{equation}\label{eqRescal}
\varphi_\lambda(x)=\lambda \Phi\left(\lambda^{(p-1)/2}|x| \right), \ \ \lambda\in (0,\infty),
\end{equation}
furnishes a family of radial steady states such that $\varphi_\lambda(0)=\lambda$. Actually, this family includes all the positive radial steady states. We point out that $\varphi_\infty$ is invariant under the above scaling.
Therefore, in light of the above, if $p>p_S$ we see that $\varphi_\lambda\to \varphi_\infty$ pointwise in $\mathbb{R}^N\setminus \{0\}$
as $\lambda\to \infty$. For completeness, let us note that for  $0<\lambda<\mu<\infty$ the following hold:
$\varphi_\lambda$ has a unique radial intersection with $\varphi_\mu$ if $p=p_S$;
$\varphi_\lambda$ has infinitely many radial intersections with $\varphi_\mu$ if $p\in (p_S,p_{JL})$;
$\varphi_\lambda<\varphi_\mu$ if $p\geq p_{JL}$.

We are now in position to state our main result.
\begin{thm}\label{thm}
If $u$ satisfies (\ref{eqAncient}) with $p\in (p_{sg},p_{JL})$ and
\begin{equation}\label{eqSandwich}
   |u(x,t)|\leq \varphi_\infty(x),\ x\in \mathbb{R}^N\setminus \{0\},\ t\leq 0,
\end{equation}
then $u\equiv 0$.
\end{thm}

For $p\in [p_S,p_{JL})\subset (p_{sg},p_{JL})$ the above theorem was proven previously by Fila and Yanagida in \cite{filaYana} by a different approach.
Roughly, they 'squeezed' $u$ between two forward self-similar solutions. We note that  forward self-similar solutions exist also for subcritical $p$, and the required properties of theirs that
are needed to show the Liouville property   are well known (see \cite[Thm. 1.1 and Cor. 1.2]{naito} and cf. also \cite[Thm. 1]{quinterAdvDE}). So,   the proof in \cite{filaYana}
applies in the above subcritical range as well. In fact, the approach in the aforementioned reference even yields the nonexistence
in the class of functions satisfying the weaker condition
\begin{equation}\label{eqAPDE}|u(x, t)| \leq (1 + \varepsilon)\varphi_\infty (x),\ x\in \mathbb{R}^N,\ t\leq 0,\end{equation} where
$\varepsilon >0$ is given by the properties of the forward self-similar solutions. Loosely speaking, time can be considered as the 'squeezing parameter' in their proof.

In contrast, our proof   does not make use of (time dependent) similarity variables.
Instead of using time dependent solutions as barriers,  we will plainly use $\varphi_\infty$ after  appropriately  'covering' its singularity   with a  piece of $\varphi_\lambda$ from (\ref{eqRescal}) (a 'surgery' type of argument in some sense). The result is a weak supersolution of (\ref{eqAncient}). Our 'squeezing parameter' will plainly  be $\lambda>0$  through the use of Serrin's sweeping principle (see \cite[Thm. 9]{mcnab} and \cite[Thm. 2.7.1]{serg} for the elliptic case) in the spirit of the sliding method \cite{bcn}. However, to be able to start such a continuity argument,
we need that $u$ is bounded. Thankfully, as it turns out, this  can be assumed without loss of generality  in light of the scaling and doubling arguments of \cite{pqs}.
On the other hand, most likely, our approach cannot be used to prove the nonexistence in the case of the weaker condition
(\ref{eqAPDE}).

Remarkably, it was shown in \cite{filaYana} that the equation in (\ref{eqAncient}) admits  positive, entire solutions of homoclinic
and heteroclinic type for $p\in (p_S, p_L)$ and $p\in [p_S,p_{JL})$, respectively, where $p_L>p_{JL}$ stands for Lepin's exponent.

In the case of positive, entire solutions of (\ref{eqAncient}), it was shown in \cite{veronGracie} that the Liouville property holds without the restriction (\ref{eqSandwich})
if $1<p<N(N + 2)/(N - 1)^2$ for $N\geq 2$ and $p>1$ for $N=1$. We point out that $p_{sg}<N(N + 2)/(N - 1)^2<p_S$ for $N\geq 3$. Actually, it was believed that the aforementioned result should
hold for all  $p\in (1,p_S)$ in analogy with the steady state problem. At first, this was shown to be true in the case of radial solutions in \cite{pqGracie}. The problem in its full generality was settled only  just recently in \cite{quinteras}.

Our approach, being elliptic in nature, carries over with only minor modifications  to establish  the following elliptic counterpart of Theorem \ref{thm} (as we will point out, the solutions in the latter can be extended for $t\in \mathbb{R}$). In contrast, the  approach of \cite{filaYana} is intrinsically parabolic and seems to be inapplicable for this purpose.
\begin{thm}\label{thm2}
If $u$ satisfies
\begin{equation}\label{eqNM}
  \Delta u+|u|^{p-1}u=0,\ z=(x,y)\in \mathbb{R}^{N+M}, \
\textrm{with}
\ p\in \left(p_{sg}(N),p_{JL}(N)\right),\ N\geq 3,\ M\geq 0,
\end{equation}
and
\begin{equation}\label{eqscoop}
|u(x,y)|\leq \varphi_\infty(x),\ x\in \mathbb{R}^N\setminus \{0\},\ y\in \mathbb{R}^M,
\end{equation}
then $u\equiv 0$.
\end{thm}

If $M=0$, the assumption (\ref{eqscoop}) implies that $u$ is a stable solution of (\ref{eqNM}) in $\mathbb{R}^N\setminus\{0\}$\footnote{We were informed of this property by L. Dupaigne  after the first version of the paper, we  borrow his argument.} (see \cite[Ch. 1]{dupaBook} for the definition).
Indeed, as in \cite[Prop. 1.3.2]{dupaBook}, it is easy to check that the difference
$\varphi_\infty-|u|
$
is a positive weak supersolution of the linearized operator
$-\Delta-p|u|^{p-1}$ in $\mathbb{R}^N\setminus \{0\}$ (see \cite[Ch. 9]{Evans} for the definition, and also \cite[Lem. I.1]{berest}).
Consequently, by the obvious weak version of \cite[Prop. 1.2.1]{dupaBook} (see also \cite[App. B]{Famu}), we infer that $u$ is stable in  $\mathbb{R}^N\setminus \{0\}$.
 Therefore, in the special case $M=0$, our result follows  from \cite[Thm. 2]{farg} which asserts that in that case (\ref{eqNM}) cannot have a nontrivial solution that is stable outside  a compact set. On the other hand, we note that this viewpoint cannot be applied for general $M>0$ because the
exponent $p_{JL}(K)$ is decreasing with respect to $K$.

An analogous Liouville type result to Theorem \ref{thm} for $p\geq p_{JL}$, which takes into account that $\varphi_\lambda<\varphi_\infty$, $\lambda \in (0,\infty)$, can be found in our recent paper
\cite{sour}. In the aforementioned work we have extended, again with a sweeping argument, the Liouville type result of Polacik and Yanagida from \cite{py} who relied on (time dependent) similarity variables and invariant manifold ideas. A version of Theorem \ref{thm2} for $p\geq p_{JL}(N)$ is contained in an extended  remark in the same paper of ours.

In the case of the critical Sobolev exponent $p=p_S$, a famous result of Caffarelli, Gidas and Spruck \cite{CGS12321} asserts that all positive solutions of the steady state problem in $\mathbb{R}^N\setminus\{0\}$
are radial (whether they have a removable singularity at the origin or not).
Using this information, Schoen \cite{Schoen}
observed that all such   solutions with a singularity at the origin can be completely classified by standard ODE phase-plane analysis. They are of the form
\begin{equation}\label{eqLogaras}
u(x)= |x|^{-\frac{N-2}{2}} v\left(\ln |x|\right),
\end{equation}
where $v$ is a positive periodic solution of
\begin{equation}\label{eqODEAI}
-v''+
\frac{(N- 2)^2}{4}
v - v
^\frac{N+2}{N-2} = 0\  \textrm{in}\  \mathbb{R}.
\end{equation}
Besides of the constant  solution $\left(\frac{N-2}{2}\right)^\frac{N-2}{2}$, which gives rise to the self-similar singular solution
$\varphi_\infty$, there is a   family of periodic solutions that can be uniquely parametrized,
up to translations, by their minimal value which spans the interval $\left(0,\left(\frac{N-2}{2}\right)^\frac{N-2}{2} \right)$.
These periodic
solutions have a unique local maximum and minimum per period. In fact, they are symmetric
with respect to their local extrema. The  singular solutions of (\ref{eqAncient})  that they produce via (\ref{eqLogaras})
are frequently called of Delaunay-type in comparison with Delaunay surfaces which are singly periodic, rotationally symmetric surfaces with constant
mean curvature (we refer to \cite{pacardiaresy} and the references therein for this connection).
We point out that each Delaunay-type singular solution has infinitely many radial intersections with $\varphi_\infty$.
Actually, the radial regular steady state $\Phi$ of (\ref{eqAncient}) is given by (\ref{eqLogaras}) with $v$ an appropriate translation
of the positive, even homoclinic solution of (\ref{eqODEAI}). Remarkably, the latter solution can be computed  explicitly and is equal to
$\left(N(N-2)\right)^\frac{N-2}{4}\left(2 \cosh(\cdot)\right)^{-\frac{N-2}{2}}$.
 Let us note in passing that the translation invariance of (\ref{eqODEAI}) echoes the scaling invariance (\ref{eqRescal}) of (\ref{eqAncient}).
It is worth mentioning that an analogous transformation to (\ref{eqLogaras}) also applies for $p\neq p_S$.
However, the corresponding
second order autonomous  ODE for $v$   is dissipative and thus has no nonconstant periodic or homoclinic solutions (it has, however, heteroclinic solutions for $p\in (p_{sg},p_S)$ that give rise to fast decaying singular solutions, see Remark \ref{remChenia} and the references therein).

Armed with the above information and by suitably adapting our approach, we can complement our main results with the following.

\begin{thm}\label{thm3Gokor}
If $p=p_S(N)$, the assertions of Theorems \ref{thm} and \ref{thm2} hold
with the righthand side of (\ref{eqSandwich}) and (\ref{eqscoop}), respectively, being an arbitrary Delaunay-type singular solution.
\end{thm}

To illustrate the delicacy of our result, at least in the parabolic case, we remark that the previously mentioned heteroclinic
solutions of \cite{filaYana} connect  $\varphi_\lambda$, $\lambda\in (0,\infty)$, as $t\to -\infty$ to the trivial solution as $t\to +\infty$
and are decreasing in time.
In the critical case $p=p_S$, it was speculated in \cite{delGrim}  that, for any
$k\in \mathbb{N}$, there exists an ancient solution that  roughly behaves like $\sum_{i=1}^{k}\varphi_{\lambda_i(t)}$ as $t\to -\infty$  for some $\lambda_i(t)\to \infty$ such that $\lambda_i/\lambda_{i+1}\to 0$.
 Clearly, such a solution must intersect  any Delaunay-type singular solution, provided that $t<0$ is sufficiently large. So, there is no contradiction with our result.

 Analogous Delaunay-type singular solutions have been studied recently in the case of the fractional Laplacian operator and for the bilaplacian one in \cite{maria} and \cite{frankie}, respectively. Whether our results can be extended in this setting is left as an interesting open problem.

The rest of the paper is essentially devoted to the proofs of our main results in the next section. In Subsection \ref{subsecGracie1}, we will prove Theorem \ref{thm}. After its proof, in Remark \ref{remToti}, we will hint at a perhaps unexpected connection  between our supersolution and a well known argument from the theory of minimal surfaces. As we have already mentioned, the proof of Theorem \ref{thm2} requires only minor modifications and  will therefore be omitted. In Subsection \ref{subsecGracie2}, we will prove Theorem \ref{thm3Gokor}. Subsequently, in Remark \ref{remChenia}, we will give a partial analog of this theorem for subcritical exponents.
Lastly, for the reader's convenience, in Appendix \ref{App} we will state a reduced version of the doubling lemma from \cite{pqs0} that is needed for our results.

\section{Proofs of the main results}In this section we will prove Theorems \ref{thm} and \ref{thm3Gokor}. In order to avoid confusion, we mention again that the proof of Theorem \ref{thm2} will be omitted as it
requires only minor adaptations.\subsection{Proof of Theorem \ref{thm}}\label{subsecGracie1}

\begin{proof}
The main idea of the proof is to construct a family of weak supersolutions of (\ref{eqAncient}) by appropriately modifying the singular solution $\varphi_\infty$ around the origin.
Our construction will hinge on the fact that, as we have already mentioned,  the radial regular steady state $\Phi$ intersects at least once with $\varphi_\infty$ since $p\in (p_{sg},p_{JL})$.
We denote by $r_1>0$ the smallest radius at which such an intersection takes place, and define a function $Z:\mathbb{R}^N\to \mathbb{R}$ with radial profile given by
\begin{equation}\label{eqZ}
Z(r)=\left\{ \begin{array}{ll}
                    \Phi(r), & 0\leq r\leq r_1, \\
                      &   \\
                    \varphi_\infty(r), & r>r_1.
                  \end{array}
\right.
\end{equation} Clearly, $Z$ is continuous  by our choice of $r_1$.
The point is that it is a weak   supersolution of (\ref{eqAncient}) (see for instance \cite[Ch. 5]{L} for the definition) because \begin{equation}\label{eq5}
                                                                             \Phi'(r_1)>\varphi_\infty'(r_1)
                                                                           \end{equation} holds (see also \cite[Lem. I.1]{berest}).
Next, according to (\ref{eqRescal}), we let
\[
z_\lambda(x)=\lambda Z\left(\lambda^{(p-1)/2}|x| \right)=\left\{ \begin{array}{ll}
                    \varphi_\lambda(r), & 0\leq r\leq s_\lambda, \\
                      &   \\
                    \varphi_\infty(r), & r>s_\lambda,
                  \end{array}
\right.\  r=|x|, \ \lambda>0,
\]
where we have denoted
\[
s_\lambda=r_1\lambda^{- (p-1)/2}.
\]
 We emphasize that we have used  that $\varphi_\infty$ is invariant under the above rescaling.
We point out that $z_\lambda\to \infty$ uniformly on $|x|\leq s_\lambda$ as $\lambda\to \infty$. On the other hand, $z_\lambda\to 0$ as $\lambda\to 0$, uniformly in $\mathbb{R}^N$.
Clearly, $z_\lambda$ is still a weak supersolution to (\ref{eqAncient}).
Actually, we will not use any   weak form of the maximum principle in the sequel (i.e. as that in \cite[Ch. 9]{Evans}). Nevertheless, the fact that $z_\lambda$ is a weak supersolution of (\ref{eqAncient})  will serve as an important guideline.

By  making partial use of our supersolution, we will first show that $u$ can be extended as a solution of the equation in (\ref{eqAncient}) for $t\in \mathbb{R}$. To this end,  the standard existence and uniqueness theory for the corresponding Cauchy problem (it is well-posed in $L^\infty(\mathbb{R}^N)$, see \cite[Prop. 51.40]{qs}) guarantees that  $u$ can be extended  in a maximal time interval of the form $(-\infty,T)$ for some $T\in (0,\infty]$.
Moreover, by the strong maximum principle for linear parabolic equations \cite[Ch. II]{L}, we assert from (\ref{eqSandwich}) that \begin{equation}\label{eqlat}
                                                                                                    |u|<\varphi_\infty, \ x\in \mathbb{R}^N\setminus \{0\},\ t\in (-\infty,T).
                                                                                                  \end{equation}
Since $u(\cdot,0)$ is a bounded function, there exists a $\lambda^*\gg 1$ such that
\[
u(x,0)<\varphi_{\lambda^*}(x),\ \ |x|\leq s_{\lambda^*}.
\]
Let $\varepsilon\in (0,T)$ be arbitrary. By virtue of the above two relations, since $u$ and $\varphi_{\lambda^*}$ are bounded on $\{|x|\leq s_{\lambda^*},\ t\in [0,T-\varepsilon]\}$,  the parabolic maximum principle \cite[Lem. 2.3]{L} (applied to the linear equation for the difference of these two solutions of (\ref{eqAncient})) yields
\[
u\leq \varphi_{\lambda^*}\ \ \textrm{for}\ |x|\leq s_{\lambda^*},\ t\in [0,T-\varepsilon].
\]
Since $\varepsilon>0$ is arbitrary, we obtain $u\leq \varphi_{\lambda^*}$ for $|x|\leq s_{\lambda^*}$, $t\in [0,T)$. Applying the same argument with $-u$ in place of $u$,
and keeping in mind (\ref{eqlat}), we conclude that $u$ remains bounded as $t\to T^-$. This means that $T=\infty$ as desired (if not, then $u$ could be continued further as a solution in contradiction to  the maximality of $T$).

Having disposed of this preliminary step, we can now turn our attention to  the Liouville property.
By nowadays standard doubling and scaling arguments \cite{pqs}, we can assume that $u$ is bounded.
In fact, we can do better and assume that\begin{equation}\label{eqbeter}
                                           |u|\leq 1\ \textrm{in}\ \mathbb{R}^N\times \mathbb{R}.
                                         \end{equation}
Indeed, let us suppose  that $|u(x_0,t_0)|>1$ for some $(x_0,t_0)\in \mathbb{R}^N\times \mathbb{R}$. Motivated from \cite{quig}, we will
 apply Lemma \ref{lemDoubling} from Appendix \ref{App} with $X=\mathbb{R}^N\times \mathbb{R}$, equipped with the parabolic distance
 \[
 d\left((x,t),(\tilde{x},\tilde{t}) \right)=|x-\tilde{x}|+\sqrt{|t-\tilde{t}|},
 \]
 and
 \[
 M(x,t)=|u|^{(p-1)/2}(x,t).
 \]
For $\texttt{y}=(x_0,t_0)$ and any $k\in \mathbb{N}$, the aforementioned lemma
   provides   $(x_k,t_k)$ such that
  \[
  M_k:=|u|^{(p-1)/2}(x_k,t_k)\geq |u|^{(p-1)/2}(x_0,t_0)
  \]
  and
  \[
  |u|^{(p-1)/2}(x,t)\leq 2 M_k\ \textrm{whenever}\ |x-x_k|+\sqrt{|t-t_k|}\leq \frac{k}{M_k}.
  \]
We note that (\ref{eqSandwich}) and the definition of $M_k$ force
\[
  M_k|x_k|\leq L^{ {(p-1)/2}}.
\]
 Hence, passing to a subsequence if necessary, we may assume that
\begin{equation}\label{eqLocalization}
M_kx_k\to y_\infty\ \textrm{for some}\ y_\infty \in \mathbb{R}^N.
\end{equation}
The rescaled functions
\[
v_k(y,s)=\rho_k^{2/(p-1)}u(x_k+\rho_ky,t_k+\rho_k^2s), \ \textrm{where}\ \rho_k=\frac{1}{2M_k},
\]
are  entire solutions of (\ref{eqAncient}) and satisfy $|v_k(0,0)|=2^{-2/(p-1)}$, $|v_k(y,s)|\leq 1$ for $|y|+\sqrt{|s|}\leq 2k$.
The parabolic regularity theory \cite[Chs. IV, VII]{L} guarantees that the sequence $\{v_k \}$ is relatively
compact in $C^{2+\theta,1+\theta/2}_{loc}$ for some $\theta\in (0,1)$. Hence, using the usual diagonal argument, passing to a further subsequence if needed, we may assume that \[v_k\to V\ \textrm{in} \ C^{2,1}_{loc}(\mathbb{R}^N\times\mathbb{R}),\] where $V$ is an entire solution to  (\ref{eqAncient})
such that  $|V|\leq 1$ and $V(0,0)\neq 0$. Furthermore, on account of (\ref{eqSandwich}), we have
\[
|v_k(y,s)|\leq \frac{L\rho_k^{2/(p-1)}}{|x_k+\rho_ky|^{2/(p-1)}}=\frac{L}{|x_k/\rho_k+y|^{2/(p-1)}}=\frac{L}{|2M_kx_k+y|^{2/(p-1)}},\ y\neq -\frac{x_k}{\rho_k}.
\]
Thus, by letting $k\to \infty$ and using  (\ref{eqLocalization}), we obtain
\[
|V(y,s)|\leq \frac{L}{|2y_\infty+y|^{2/(p-1)}}, \ y\neq -2y_\infty.
\]
Now, the spatially shifted solution
\[
W(y,s)=V(y-2y_\infty,s)
\]
satisfies $|W|\leq 1$,   $W(2y_\infty,0)\neq 0$ and (\ref{eqSandwich}).
Consequently, it is sufficient to prove the theorem for entire solutions that satisfy (\ref{eqSandwich}) with $t\in \mathbb{R}$ and (\ref{eqbeter}).
This task will take up the rest of the proof.

The main tool in the proof is Serrin's sweeping principle (see \cite[Thm. 9]{mcnab} and  \cite[Thm. 2.7.1]{serg} for the elliptic case) using the family of supersolutions $\{z_\lambda\}$.
Since $u$ is bounded and satisfies (\ref{eqSandwich}), there exists a    $\bar{\lambda} \gg 1$ such that
\[u\leq z_{\mu},\ x\in \mathbb{R}^N, \ t\in \mathbb{R},\ \textrm{for any}\ \mu\geq \bar{\lambda}.\]
Starting from $\bar{\lambda}$, we  proceed to decrease $\lambda$ while keeping the above ordering. There are only two possibilities.
Either we  can  continue all the way until we reach $\lambda=0$ or we will get 'stuck' at some first $\lambda_0>0$ and cannot continue further.
Our goal is to show that the latter scenario (to be described in more detail below) cannot happen. This will imply that $u\leq 0$. Then, the assertion of the theorem follows readily by carrying out the same procedure with $-u$ in place of $u$.

Let us suppose, to the contrary, that there exists some $\lambda_0\in (0,\bar{\lambda}]$ where we get stuck in the sense that
the set
\[
  \Lambda=\left\{\lambda\geq 0\ :\ z_\mu \geq u\ \textrm{in}\ \mathbb{R}^N\times \mathbb{R} \ \textrm{for every}\ \mu \geq \lambda\right\}
\]
coincides with $[\lambda_0,\infty)$ (by its definition $\Lambda$ is a semi infinite interval, while it is closed thanks to the continuity of $z_\mu$ with respect to $\mu$).
Clearly,   we have
\begin{equation}\label{eqcontra1}
u\leq z_{\lambda_0},\ x\in \mathbb{R}^N, \ t\in  \mathbb{R}.
\end{equation}
Keeping in mind that $z_\lambda$ depends nontrivially on $\lambda$ only in the space-time cylinder $\{ |x|<s_\lambda,\ t\in \mathbb{R}\}$ (where it is equal to $\varphi_\lambda$), and (\ref{eqlat}) with $T=\infty$, we get $\lambda_k\in (0,\lambda_0)$ such that  $\lambda_k\to {\lambda_0}$ as $k\to \infty$,    $x_k\in\mathbb{R}^N$ with $|x_k|< s_{\lambda_k}$, and $t_k\in \mathbb{R}$ such that
  \begin{equation}\label{eqcontra2}
    u(x_k,t_k)>\varphi_{\lambda_k}(x_k),
  \end{equation}
(the reader should not  be confused with the repeated use of notation in different contexts within the proof).
The whole argument is actually reminiscent to the famous sliding method \cite{bcn} for elliptic problems, when translating  a compactly supported subsolution (as in \cite[Thm. 2.1]{danceria} for instance). We also note that $z_\lambda$ and $z_\mu$ with $\lambda< \mu$ may intersect each other in $|x|<s_\mu$  as is the case in the aforementioned procedure.
Passing to a subsequence if necessary, we may  assume that
\begin{equation}\label{eqball}
  x_k\to x_\infty\  \textrm{for some}\ x_\infty\in \mathbb{R}^N\ \textrm{such that}\ |x_\infty|\leq s_{\lambda_0}.
\end{equation}

If the sequence $\{t_k\}$ is bounded, passing to a further subsequence if needed, we may assume that $t_k\to t_\infty$ for some
$t_\infty\in \mathbb{R}$. From (\ref{eqcontra1}) and (\ref{eqcontra2}), it follows that  $u(x_\infty,t_\infty)= \varphi_{\lambda_0}(x_\infty).$
In fact, thanks to (\ref{eqlat}) with $T=\infty$, we see that  $|x_\infty|\neq s_{\lambda_0}$. Thus, by virtue of (\ref{eqcontra1}) and the parabolic strong maximum principle \cite{Evans,L} (applied in the linear equation for the difference $u-\varphi_{\lambda_0}$ sufficiently close to $(x_\infty,t_\infty)$), we deduce that $u$ coincides with $\varphi_{\lambda_0}$ in some neighborhood of $(x_\infty,t_\infty)$. In turn, by repeated applications of the strong maximum principle, we obtain $u\equiv \varphi_{\lambda_0}$ which is clearly absurd on account of (\ref{eqSandwich}).

It remains to deal with the case where, up to a subsequence,   $t_k\to -\infty$ (the case where $t_k \to +\infty$ can be handled similarly). To this end, we consider the time translated solutions
\[
u_k(x,t)=u(x,t+t_k),\ x\in \mathbb{R}^N,\ t\in \mathbb{R}.
\]
From (\ref{eqcontra1}) and (\ref{eqcontra2}) it follows that \[z_{\lambda_0}\geq u_k\ \textrm{in}\ \mathbb{R}^N\times\mathbb{R}\ \textrm{and}\  u_k(x_k,0)>\varphi_{\lambda_k}(x_k),\] respectively.
Since $u$ is bounded, as before, by the usual  diagonal-compactness argument, possibly up to a  further subsequence,  we have
\[
u_k\to U\ \textrm{in}\ C^{2,1}_{loc}(\mathbb{R}^N\times \mathbb{R}),
\]
where $U$ is an entire solution to (\ref{eqAncient}) such that
\[
  U\leq z_{\lambda_0} \ \textrm{in}\ \mathbb{R}^N\times\mathbb{R}\ \textrm{and}\  U(x_\infty,0)\geq \varphi_{\lambda_0}(x_\infty).
\]
In particular, we get
$
  U(x_\infty,0)= \varphi_{\lambda_0}(x_\infty).
$ Intuitively, keeping in mind (\ref{eq5}), it is clear that we have been led to a contradiction. The rigorous justification is easy.
Indeed, by virtue of  (\ref{eq5}),  the above imply that
\[
U\leq \varphi_{\lambda_0}, \ |x|<s_{\lambda_0}+\delta,\ |t|<1,
\]
for some sufficiently small  $\delta>0$. Then, as before,   we deduce by the strong maximum principle that $U\equiv \varphi_{\lambda_0}$ which is impossible.
\end{proof}

\begin{rem}\label{remToti}
In \cite{sour},  motivated mainly from \cite{cabrouk}, we highlighted a heuristic connection of (\ref{eqAncient}) to ancient solutions of the mean curvature flow. In that context, our time independent supersolution in (\ref{eqZ}) relates to the competitor that is  used in order to show that the symmetric minimal cones are not area minimizers in low dimensions (see for instance \cite{angi}). It is worth noting that, in analogy to our proof, a  sweeping argument with the aforementioned competitor is also possible as in the alternative proof of \cite[Thm. 1.8]{cabrouk}.
\end{rem}\subsection{Proof of Theorem \ref{thm3Gokor}}\label{subsecGracie2}

\begin{proof}The proof is similar to that of Theorems \ref{thm} and \ref{thm2} apart from some technical modifications.
We will give a sketch of the proof only for the parabolic problem (the elliptic case is analogous) and point out the main differences.

Let us denote by
\[\psi(r)=h(\ln r)r^{-\frac{N-2}{2}},\ r=|x|> 0,\]
with $h>0$ a $T$-periodic solution of (\ref{eqODEAI}),
 a Delaunay-type singular solution that bounds the absolute value of $u$.
For each $\lambda\in (0,\infty)$, the homoclinic solution of (\ref{eqODEAI}) that gives $\varphi_\lambda$ via (\ref{eqLogaras}) intersects at least twice with $h$ (this can be seen easily from the phase plane portrait). Hence, there exists a first radius $\tau_\lambda>0$
at which $\varphi_\lambda$ and $\psi$ intersect.  Clearly, we have $\tau_\lambda\to 0$ as $\lambda \to \infty$ and $\tau_\lambda\to \infty$ as $\lambda \to 0$. Moreover, since such an intersection is transverse (by the uniqueness of the IVP for the radial ODE),
the implicit function theorem guarantees that $\tau_\lambda$ varies smoothly with respect to $\lambda>0$.
Keeping in mind that $\psi$ is not invariant under the scaling in (\ref{eqRescal}) (unless $h\equiv L$ of course), we now define our supersolution $z_\lambda$ directly as
\[
z_\lambda(x)=\left\{\begin{array}{ll}\varphi_\lambda(r),& 0\leq r\leq \tau_\lambda,\\
&\\
\psi(r),& r>\tau_\lambda,
\end{array} \right. \ r=|x|.
\]
As before, we can use $z_\lambda$ as a barrier in order to show that $u$ cannot blow up in finite time. Therefore, we may again assume that $u$ is an entire solution
to (\ref{eqAncient}) that satisfies   \begin{equation}\label{eqDivina} |u(x,t)|\leq h(\ln |x|)|x|^{-\frac{N-2}{2}},\ \ x\in \mathbb{R}^N\setminus\{0\},\ t\in \mathbb{R}.\end{equation}
As in the proof of Theorem \ref{thm},  by applying Serrin's sweeping principle,  we can conclude that $u\equiv 0$ under the additional assumption that it is bounded.

It  remains to verify that, by the doubling lemma as in the proof of Theorem \ref{thm}, we can assume without loss of generality  that (\ref{eqbeter}) holds.
To this end, assuming that this was not the case, we define $M$, $(x_k,t_k)$, $M_k$, $\rho_k$ and $v_k(y,s)$ analogously to the aforementioned proof. We quickly come across a  minor difference which is that now we have
\[
M_k|x_k|\leq \|h\|_{L^\infty(\mathbb{R})}^\frac{2}{N-2}.
\]
Nevertheless, up to a subsequence, we still have $M_kx_k\to y_\infty$ for some $y_\infty \in \mathbb{R}^N$.
Moreover, we still have  the local convergence  of $v_k$ to some bounded, nontrivial limiting solution $V$.
However, the main differences arise when passing to the limit in the rescaled form of  (\ref{eqDivina}).
More precisely,  the latter gives
\[
|v_k(y,s)|\leq \frac{\rho_k^\frac{N-2}{2}h\left(\ln |x_k+\rho_k y|\right)}{|x_k+\rho_k y|^{\frac{N-2}{2}}}=
\frac{h\left(\ln |x_k+\rho_k y|\right)}{|x_k/\rho_k+ y|^{\frac{N-2}{2}}}
\]
for $y\neq -x_k/\rho_k$, $s\in \mathbb{R}$. Based on the identity
\[
\ln |x_k+\rho_k y|=\ln\rho_k+ \ln|\frac{x_k}{\rho_k}+y|,
\]
we   decompose $\ln\rho_k$ as
\[
\ln\rho_k=m_k T+d_k,
\]
with $m_k \in \mathbb{Z}$ and $|d_k|\leq T$. By passing to a further subsequence if necessary, we may assume that $d_k\to d_\infty$ for some $d_\infty \in \mathbb{R}$.
Since $h$ is $T$-periodic, we obtain
\[h\left(\ln |x_k+\rho_k y|\right)=
h\left(m_k T+d_k+ \ln|\frac{x_k}{\rho_k}+y|\right)=h\left(d_k+ \ln|\frac{x_k}{\rho_k}+y|\right).
\]
Consequently, recalling the definition of $M_k$, we get
\[
|v_k(y,s)|\leq
\frac{h\left(d_k+ \ln|2M_kx_k+y|\right)}{|2M_kx_k+ y|^{\frac{N-2}{2}}}, \ \ y\neq -2M_kx_k.
\]
Letting $k\to \infty$, we deduce that
\[
|V(y,s)|\leq
\frac{h\left(d_\infty+ \ln|2y_\infty+y|\right)}{|2y_\infty+ y|^{\frac{N-2}{2}}}=\frac{h\left( \ln\left(e^{d_\infty}|2y_\infty+y|\right)\right)}{|2y_\infty+ y|^{\frac{N-2}{2}}}, \ \ y\neq -2y_\infty.
\]
We remark that the righthand side of the above relation is plainly a rescaling (according to (\ref{eqRescal})) and a translation of the Delaunay-type singular solution $\psi$.
In other words, after a translation, $V$ satisfies (\ref{eqDivina}) with $h$ replaced by  a positive  $e^{-d_\infty}T$-periodic solution of (\ref{eqODEAI}).
Hence, $V$ is a bounded solution that satisfies the assumptions of the theorem, which is what we wanted.
\end{proof}
\begin{rem}\label{remChenia}
If $p\in(p_{sg},p_S)$, for any $a>0$, there exists a positive, radial singular solution $\phi$ to the steady state problem such that
\[
\lim_{r\to 0}r^{\frac{2}{p-1}}\phi(r)=L\ \textrm{and}\ \lim_{r\to \infty}r^{N-2}\phi(r)=a
\]
(see \cite[Prop. 2.2]{cheniaresy}). We note that these singular solutions decay faster than the self-similar one
as $|x|\to \infty$. In fact, the aforementioned ones are the only positive solutions of the steady state problem in
$\mathbb{R}^N\setminus\{0\}$ with $p$ in this range (see \cite[Prop. 3.1]{seriniaresy}).
We observe that $\varphi_\lambda$ with $\lambda \in (0,\infty)$ must intersect at least twice with each such fast decaying singular solution.
Indeed, if not then by the discussion following Theorem \ref{thm2} we would have that $\varphi_\lambda$ is a stable solution of the steady state problem in its support which is absurd
(see for instance \cite[Ex. 1.2.3]{dupaBook}). In light of this property, by arguing as in the proof of Theorem \ref{thm3Gokor} we can show  that the Liouville property holds
for bounded, ancient solutions to (\ref{eqAncient}) that are smaller in absolute value than such a fast decaying singular steady state.
\end{rem}

\appendix
\section{A doubling lemma from \cite{pqs0}}\label{App}
In this small appendix, we will state for the reader's convenience the following reduced version of \cite[Lem. 5.1]{pqs0} that we referred to in the proof of Theorem \ref{thm}.

\begin{lem}\label{lemDoubling}
Let $(X,d)$ be a complete metric space and
 $M :  X\to [0,\infty)$ be bounded on compact subsets of $X$. Fix a $\texttt{y}\in X$ such that $M(\texttt{y})>0$ and a real $k>0$. Then, there exists $\texttt{x}\in X$ such that
\[
   M(\texttt{x})\geq M(\texttt{y})\]
and
\[
  M(\texttt{z})\leq 2 M(\texttt{x})\ \ \textrm{whenever}\ \ d(\texttt{z},\texttt{x})\leq
   \frac{k}{M(\texttt{x})}.
\]
\end{lem}
\begin{rem}
Our formulation of the doubling lemma is restricted to the whole metric space $X$ following a related comment   in \cite[Sec. 2]{quig}.
We also note that we assume $M$ to be nonnegative instead of strictly positive, as was the case in the previous references. However,  if $M(\texttt{y})>0$ then throughout the  proof of \cite[Lem. 5.1]{pqs0} we observed that $M$ is evaluated only at points where $M\geq M(\texttt{y})$. Thus, there is no loss of generality. In fact, the doubling lemma, as stated in the aforementioned reference, has been previously applied with $M\geq 0$ (possibly taking zero values)
in \cite[Thm. 1.7]{bpq} and \cite[Lem. 5.1]{dav}.
\end{rem}

\textbf{Acknowledgments.} The author would like to thank IACM of FORTH, where this paper was written, for the hospitality. This work has received funding from the Hellenic Foundation for Research and Innovation (HFRI) and the General Secretariat for Research and Technology (GSRT), under grant agreement No 1889.

\end{document}